\documentclass[11pt]{amsart}%
\usepackage{graphicx,amscd,color,amsmath,amsfonts,amssymb,geometry}
\usepackage{amsmath}
\usepackage{amsfonts}
\usepackage{amssymb}
\usepackage{graphicx}%
\providecommand{\U}[1]{\protect\rule{.1in}{.1in}}
\newtheorem{theorem}{Theorem}[section]

\begin{document}
\title[Optimal Hardy--Littlewood constants for $2$-homogeneous polynomials on
$\ell_{p}(\mathbb{R}^{2})$ ]{The optimal Hardy--Littlewood constants for $2$-homogeneous polynomials on
$\ell_{p}(\mathbb{R}^{2})$ for $2<p<4$ are $2^{2/p}$}
\author[W. Cavalcante]{W. Cavalcante}
\address{Departamento de Matem\'{a}tica, \\
\indent Universidade Federal de Pernambuco, \\
\indent50.740-560 - Recife, Brazil.}
\email{wasthenny.wc@gmail.com}
\author[D. N\'{u}\~{n}ez]{D. N\'{u}\~{n}ez-Alarc\'{o}n}
\address{Departamento de Matem\'{a}tica, \\
\indent Universidade Federal de Pernambuco, \\
\indent50.740-560 - Recife, Brazil.}
\email{danielnunezal@gmail.com}
\author[D. Pellegrino]{D. Pellegrino}
\address{Departamento de Matem\'{a}tica, \\
\indent Universidade Federal da Para\'{\i}ba, \\
\indent 58.051-900 - Jo\~{a}o Pessoa, Brazil.}
\email{pellegrino@pq.cnpq.br and dmpellegrino@gmail.com}
\thanks{W. Cavalcante was supported by Capes D. N\'{u}\~{n}ez was supported by CNPq
Grant 461797/2014-3. D. Pellegrino was supported by CNPq Grant 477124/2012-7
and INCT-Matem\'{a}tica.}
\keywords{Optimal constants; Hardy--Littlewood inequality}
\subjclass[2010]{11Y60, 47B10, 26D15, 46B25.}

\begin{abstract}
We show that the optimal constants for the Hardy--Littlewood inequalities for
$2$-homogeneous polynomials on $\ell_{p}(\mathbb{R}^{2})$ are precisely
$2^{2/p}$ for all $2<p<4.$

\end{abstract}
\maketitle

%\thanks{2010 Mathematics Subject Classification:}

\section{Introduction}

The Hardy--Littlewood inequality for bilinear forms in $\ell_{p}$ spaces were
proved in 1934 \cite{hardy}; these inequalities and the classical
Bohnenblust--Hille inequality \cite{bh} consist in optimal extensions of
Littlewood's $4/3$ inequality \cite{LLL} (originally stated for $c_{0}$
spaces). In the last years the interest in this subject (which can be
considered part of the theory of multiple summing and absolutely summing
operators) was renewed with applications in several fields of Mathematics and
even in Physics (see \cite{bps, ann, monta}) and several authors became
interested in this field (\cite{abps, apd, cara, mic, def33, sur, dimant, gal, pell,
popa, pra, diana}). Very recently the subject became to be investigated via
numerical and computational assistance due to several challenging problems
which seemed quite difficult to be solved analytically (without computer assistance).\

There is no doubt that computational assistance is important in this subject
but some problems arise with this approach: for instance, the concrete
estimates of the constants obtained with computational assistance are just
approximations (even if we have thousands of decimal digits of confidence) of
the exact values of the constants and closed (and elegant) formulas for the
optimal constants are difficult (or even essentially impossible) to be
achieved just with the help of computers.

As mentioned before, the search of optimal constants for these kind of
inequalities has important applications but, as a matter of fact, even for
$2$-homogeneous polynomials the optimal constants are unknown. The same
happens with the constants of the Hardy--Littlewood inequalities. In this note
we obtain simple formulas for the optimal Hardy--Littlewood constants for
$2$-homogeneous polynomials on $\ell_{p}(\mathbb{R}^{2})$ for all $2<p\leq4.$
Up to now, the only known simple (explicit) formula for these constants is
$2^{1/2}$ for $p=4$, due to Araujo \textit{et al. \cite{aaww}, }extending
previous results of \cite{waaa}.

For $\mathbb{K}$ be $\mathbb{R}$ or $\mathbb{C}$ and $\alpha=(\alpha
_{1},\ldots,\alpha_{n})\in{\mathbb{N}}^{n}$, we define $|\alpha|:=\alpha
_{1}+\cdots+\alpha_{n}$. By $\mathbf{x}^{\alpha}$ we shall mean the monomial
$x_{1}^{\alpha_{1}}\cdots x_{n}^{\alpha_{n}}$ for any $\mathbf{x}%
=(x_{1},\ldots,x_{n})\in{\mathbb{K}}^{n}$. The polynomial version of
Littlewood's $4/3$ theorem asserts that, given $n\geq1$, there is a constant
$B_{\mathbb{K},2}^{\mathrm{pol}}\geq1$ such that
\[
\left(  {\sum\limits_{\left\vert \alpha\right\vert =2}}\left\vert a_{\alpha
}\right\vert ^{\frac{4}{3}}\right)  ^{\frac{3}{4}}\leq B_{\mathbb{K}%
,2}^{\mathrm{pol}}\left\Vert P\right\Vert
\]
for all $2$-homogeneous polynomials $P:$ $\ell_{\infty}^{n}\rightarrow
\mathbb{K}$ given by
\[
P(x_{1},...,x_{n})=\sum_{|\alpha|=2}a_{\alpha}\mathbf{{x}^{\alpha},}%
\]
and all positive integers $n$, where $\Vert P\Vert:=\sup_{z\in B_{\ell
_{\infty}^{n}}}|P(z)|$. It is well-known that the exponent $\frac{4}{3}$ is sharp.

The change of $\ell_{\infty}^{n}$ by $\ell_{p}^{n}$ gives us the polynomial
Hardy--Littlewood inequality whose optimal exponents are $\frac{4p}{3p-4}$ for
$4\leq p\leq\infty$ and $\frac{p}{p-2}$ for $2<p\leq4$. More precisely, given
$n\geq1$, there is a constant $C_{\mathbb{K},2,p}^{\mathrm{pol}}\geq1$ such
that
\[
\left(  {\sum\limits_{\left\vert \alpha\right\vert =2}}\left\vert a_{\alpha
}\right\vert ^{\frac{4p}{3p-4}}\right)  ^{\frac{3p-4}{4p}}\leq C_{\mathbb{K}%
,2,p}^{\mathrm{pol}}\left\Vert P\right\Vert ,
\]
for all $2$-homogeneous polynomials on $\ell_{p}^{n}$ with $4\leq p\leq\infty$
given by $P(x_{1},\ldots,x_{n})=\sum_{|\alpha|=m}a_{\alpha}\mathbf{{x}%
^{\alpha}}$. When $2<p\leq4$ we have%
\[
\left(  {\sum\limits_{\left\vert \alpha\right\vert =2}}\left\vert a_{\alpha
}\right\vert ^{\frac{p}{p-2}}\right)  ^{\frac{p-2}{p}}\leq C_{\mathbb{K}%
,2,p}^{\mathrm{pol}}\left\Vert P\right\Vert .
\]
The main result of this note shows that $C_{\mathbb{R},2,p}^{\mathrm{pol}%
}=2^{2/p}$ when we are restricted to the case $n=2$ and $2<p\leq4$ (as
mentioned before, the case $p=4$ is already known). More precisely, our main
result is the following:

\begin{theorem}
\label{800}If $2<p\leq4$, then the optimal Hardy--Littlewood constants for
$2$-homogeneous polynomials on $\ell_{p}(\mathbb{R}^{2})$ are $2^{\frac{2}{p}%
}$.
\end{theorem}

\section{Proof of the main result: part 1}

\bigskip\

The following result \bigskip due to B. Grecu \cite{grecu} will be crucial for
our goals:

\begin{theorem}
\label{gre}For $p>2$, a $2$-homogeneous norm one polynomial $P$ is a extreme
point of the unit ball of $P(^{2}\ell_{p}^{2})$ if, and only if,

(i) $P(x,y)=ax^{2}+cy^{2}$, with $ac>0$ and $\left\Vert \left(  a,c\right)
\right\Vert _{\frac{p}{p-2}}=1$ or

(ii) $P(x,y)=\pm\left(  \frac{a^{p}-b^{p}}{a+b^{2}}\left(  x^{2}-y^{2}\right)
+2ab\frac{a^{p-2}+b^{p-2}}{a^{2}+b^{2}}xy\right)  ,$ with $a,b>0$ and
$\left\Vert \left(  a,b\right)  \right\Vert _{p}=1.$
\end{theorem}

\bigskip We know that for all $2$-homogeneous polynomials $P:$ $\ell_{p}%
^{n}\rightarrow\mathbb{K}$ given by
\[
P(x_{1},...,x_{n})=\sum_{|\alpha|=2}a_{\alpha}\mathbf{{x}^{\alpha},}%
\]
the formula%
\[
\left\vert P\right\vert _{q}:=\left(  \sum_{\left\vert \alpha\right\vert
=2}\left\vert a_{\alpha}\right\vert ^{q}\right)  ^{\frac{1}{q}}%
\]
defines a norm for all $q\geq1$. Since $\ell_{p}^{n}$ is finite-dimensional,
it is obvious that $\left\Vert .\right\Vert $ and $\left\vert .\right\vert
_{q}$ are equivalent; so there are constants $C_{2,n,p,q}$ such that
\begin{equation}
\left\vert P\right\vert _{q}\leq C_{2,n,p,q}\left\Vert P\right\Vert \label{bh}%
\end{equation}
for all $P\in\mathcal{P}$ $\left(  ^{2}\ell_{p}^{n}\right)  $. We shall
investigate $C_{2,n,p,q}$ in the particular case in which $n=2$, and $p>2$,
and $q\geq1$, i.e., we shall investigate $C_{2,2,p,q}.$ The following equality
(for $p>2$ and $q\geq1$)
\begin{equation}
C_{2,2,p,q}=\max_{a\in\left[  0,1\right]  }\left[  \left(  2\left\vert
\frac{2a^{p}-1}{a^{2}+\left(  1-a^{p}\right)  ^{2/p}}\right\vert ^{q}+\left(
2a\left(  1-a^{p}\right)  ^{\frac{1}{p}}\frac{a^{p-2}+\left(  1-a^{p}\right)
^{\frac{p-2}{p}}}{a^{2}+\left(  1-a^{p}\right)  ^{2/p}}\right)  ^{q}\right)
^{\frac{1}{q}}\right]  \label{7779}%
\end{equation}
due to Araujo \textit{et al}. (\cite{aaww}), is also important for our goals.
We present a proof of (\ref{7779}) for the sake of completeness.

From the Krein--Milman Theorem it is well-known that the optimal constants
$C_{2,2,p,q}$ shall be searched within extreme polynomials. So, using Theorem
\ref{gre}, \ we conclude that for $2<p$ and $q\geq1$ we have
\[
C_{2,2,p,q}=\max\left[
\begin{array}
[c]{c}%
\max_{a\in\left[  0,1\right]  }\left[  a^{q}+\left(  1-a^{\frac{p}{p-2}%
}\right)  ^{\frac{p-2}{p}q}\right]  ^{\frac{1}{q}},\\
\max_{a\in\left[  0,1\right]  }\left[  \left(  2\left\vert \frac{2a^{p}%
-1}{a^{2}+\left(  1-a^{p}\right)  ^{2/p}}\right\vert ^{q}+\left(  2a\left(
1-a^{p}\right)  ^{\frac{1}{p}}\frac{a^{p-2}+\left(  1-a^{p}\right)
^{\frac{p-2}{p}}}{a^{2}+\left(  1-a^{p}\right)  ^{2/p}}\right)  ^{q}\right)
^{\frac{1}{q}}\right]
\end{array}
\right]
\]
Note that%
\begin{equation}
\max_{a\in\left[  0,1\right]  }\left[  a^{q}+\left(  1-a^{\frac{p}{p-2}%
}\right)  ^{\frac{p-2}{p}q}\right]  ^{\frac{1}{q}}\leq2^{\frac{2}{p}%
}.\label{ip}%
\end{equation}
In fact, since $\left\Vert \cdot\right\Vert _{q}\leq\left\Vert \cdot
\right\Vert _{1},$ we have
\[
\left[  a^{q}+\left(  1-a^{\frac{p}{p-2}}\right)  ^{\frac{p-2}{p}q}\right]
^{\frac{1}{q}}\leq\left[  a+\left(  1-a^{\frac{p}{p-2}}\right)  ^{\frac
{p-2}{p}}\right]  .
\]

On the other hand, fixing $p>2$ and deriving the function $f\left(  a\right)
=a+\left(  1-a^{\frac{p}{p-2}}\right)  ^{\frac{p-2}{p}}$ we have
\[
1+-\frac{p-2}{p}\left(  1-a^{\frac{p}{p-2}}\right)  ^{\frac{-2}{p}}\left(
\frac{p}{p-2}a^{\frac{2}{p-2}}\right)  =f^{\prime}\left(  a\right)  ,
\]
Note that $f^{\prime}\left(  a\right)  $ is well-defined for all $a\in\left(
0,1\right)  \,$, and to solve $f^{\prime}\left(  a\right)  =0$ is equivalent
to solve
\[
a^{\frac{2}{p-2}}=\left(  1-a^{\frac{p}{p-2}}\right)  ^{\frac{2}{p}},
\]
and the unique real value that verifies this equality is $a_{0}=2^{\frac
{2-p}{p}}$, and
\[
1=f\left(  0\right)  =f\left(  1\right)  <f\left(  a_{0}\right)  =2^{\frac
{2}{p}}.
\]
To summarize, for each $p$, $f\left(  a\right)  $ attains its maximum in
$\left(  0,1\right)  $ at $a_{0}=2^{\frac{2-p}{p}}$. We have%
\[
\max_{a\in\left[  0,1\right]  }\left[  a^{q}+\left(  1-a^{\frac{p}{p-2}%
}\right)  ^{\frac{p-2}{p}q}\right]  ^{\frac{1}{q}}\leq\max_{a\in\left[
0,1\right]  }\left[  a+\left(  1-a^{\frac{p}{p-2}}\right)  ^{\frac{p-2}{p}%
}\right]  =2^{\frac{2}{p}},
\]
for all $q\geq1,$ and we obtain (\ref{ip}). On the other hand, estimating the
function
\[
\left(  2\left\vert \frac{2a^{p}-1}{a^{2}+\left(  1-a^{p}\right)  ^{2/p}%
}\right\vert ^{q}+\left(  2a\left(  1-a^{p}\right)  ^{\frac{1}{p}}%
\frac{a^{p-2}+\left(  1-a^{p}\right)  ^{\frac{p-2}{p}}}{a^{2}+\left(
1-a^{p}\right)  ^{2/p}}\right)  ^{q}\right)  ^{\frac{1}{q}}%
\]
at the point $a_{1}=2^{\frac{-1}{p}}$ we obtain
\begin{equation}
\left(  2\left\vert \frac{2a_{1}^{p}-1}{a_{1}^{2}+\left(  1-a_{1}^{p}\right)
^{2/p}}\right\vert ^{q}+\left(  2a_{1}\left(  1-a_{1}^{p}\right)  ^{\frac
{1}{p}}\frac{a_{1}^{p-2}+\left(  1-a_{1}^{p}\right)  ^{\frac{p-2}{p}}}%
{a_{1}^{2}+\left(  1-a_{1}^{p}\right)  ^{2/p}}\right)  ^{q}\right)  ^{\frac
{1}{q}}=2^{\frac{2}{p}}\label{pyt}%
\end{equation}
and the proof of (\ref{7779}) is done.

\section{\bigskip Proof of the main theorem: part 2}

\bigskip It suffices to prove that if $2<p\leq4$ and $q\geq2$, then
\[
\max_{a\in\left[  0,1\right]  }\left[  \left(  2\left\vert \frac{2a^{p}%
-1}{a^{2}+\left(  1-a^{p}\right)  ^{2/p}}\right\vert ^{q}+\left(  2a\left(
1-a^{p}\right)  ^{\frac{1}{p}}\frac{a^{p-2}+\left(  1-a^{p}\right)
^{\frac{p-2}{p}}}{a^{2}+\left(  1-a^{p}\right)  ^{2/p}}\right)  ^{q}\right)
^{\frac{1}{q}}\right]  =2^{\frac{2}{p}}.
\]

We first prove the case $2<p\leq4$ and $q=2$. We have just seen in (\ref{pyt})
that
\begin{align*}
&  \left(  2\left\vert \frac{2\left(  2^{-1/p}\right)  ^{p}-1}{a^{2}+\left(
1-\left(  2^{-1/p}\right)  ^{p}\right)  ^{2/p}}\right\vert ^{2}+\left(
2\left(  2^{-1/p}\right)  \left(  1-\left(  2^{-1/p}\right)  ^{p}\right)
^{\frac{1}{p}}\frac{\left(  2^{-1/p}\right)  ^{p-2}+\left(  1-\left(
2^{-1/p}\right)  ^{p}\right)  ^{\frac{p-2}{p}}}{\left(  2^{-1/p}\right)
^{2}+\left(  1-\left(  2^{-1/p}\right)  ^{p}\right)  ^{2/p}}\right)
^{2}\right)  ^{\frac{1}{2}}\\
&  =2^{\frac{2}{p}}%
\end{align*}

On the other hand, from
\[
C_{2,2,p,2}=\max_{a\in\left[  0,1\right]  }\left\{  \left(  2\left\vert
\frac{2a^{p}-1}{a^{2}+\left(  1-a^{p}\right)  ^{2/p}}\right\vert ^{2}+\left(
2a\left(  1-a^{p}\right)  ^{\frac{1}{p}}\frac{a^{p-2}+\left(  1-a^{p}\right)
^{\frac{p-2}{p}}}{a^{2}+\left(  1-a^{p}\right)  ^{2/p}}\right)  ^{2}\right)
^{\frac{1}{2}}\right\}  ,
\]
we will see that%
\[
C_{2,2,p,2}=2^{\frac{2}{p}}.
\]

In fact, defining for each $2<p\leq4$ the function $g:[0,1]\rightarrow
\mathbb{R}$ given by
\[
g\left(  a\right)  =2\left\vert \frac{2a^{p}-1}{a^{2}+\left(  1-a^{p}\right)
^{2/p}}\right\vert ^{2}+\left(  2a\left(  1-a^{p}\right)  ^{\frac{1}{p}}%
\frac{a^{p-2}+\left(  1-a^{p}\right)  ^{\frac{p-2}{p}}}{a^{2}+\left(
1-a^{p}\right)  ^{2/p}}\right)  ^{2}%
\]
and estimating its derivative we obtain%

\begin{align*}
&  g^{\prime}\left(  a\right) \\
&  =\frac{-8\left(  a^{p}\left(  1-a^{p}\right)  ^{\frac{4}{p}}+a^{p+4}%
-a^{4}\right)  \left(  a^{p}\left(  -a^{p}+1\right)  ^{\frac{2}{p}}\left(
p-a^{p}p+2a^{p}-1\right)  -\allowbreak a^{2}\left(  a^{p}-1\right)  \left(
a^{p}p-2a^{p}+1\right)  \right)  }{a^{3}\left(  \left(  1-a^{p}\right)
^{\frac{2}{p}}+a^{2}\right)  ^{3}\left(  a^{p}-1\right)  \left(
1-a^{p}\right)  ^{\frac{2}{p}}}\allowbreak.
\end{align*}

Now we observe that $g^{\prime}\left(  a\right)  $ is well-defined for all
$a\in\left(  0,1\right)  $ and moreover has precisely one zero in the interval
$\left(  0,1\right)  ,$ and this zero is attained at $2^{-1/p}$. In fact, if
$a\not =0$, $a\not =1$, $a\not =$ $2^{-1/p}$, we have $g^{\prime}\left(
a\right)  =0$ if, and only if,
\[
a^{p}\left(  -a^{p}+1\right)  ^{\frac{2}{p}}\left(  p-a^{p}p+2a^{p}-1\right)
-\allowbreak a^{2}\left(  a^{p}-1\right)  \left(  a^{p}p-2a^{p}+1\right)  =0
\]
and this never happens in $\left(  0,1\right)  $, because in this interval%
\begin{equation}
a^{p}\left(  -a^{p}+1\right)  ^{\frac{2}{p}}\left(  p-a^{p}p+2a^{p}-1\right)
-\allowbreak a^{2}\left(  a^{p}-1\right)  \left(  a^{p}p-2a^{p}+1\right)  >0.
\label{776}%
\end{equation}
The reader can be convinced of (\ref{776}) by observing that (\ref{776}) \ is
equivalent to prove that%
\[
a^{p}\left(  -a^{p}+1\right)  ^{\frac{2}{p}}\left(  p-a^{p}p+2a^{p}-1\right)
>\allowbreak a^{2}\left(  a^{p}-1\right)  \left(  a^{p}p-2a^{p}+1\right)
\]
and, this is equivalent to prove that
\[
-a^{p}\left(  -a^{p}+1\right)  ^{\frac{2-p}{p}}\left(  p-a^{p}p+2a^{p}%
-1\right)  <\allowbreak a^{2}\left(  a^{p}p-2a^{p}+1\right)  .
\]
But this last inequality is true since the left-hand-side is always negative
while the right-hand-side is always positive in $\left(  0,1\right)  $, for
$p>2.$ Thus we conclude that $g^{\prime}\left(  a\right)  $ has exactly one
zero in the interval $\left(  0,1\right)  ,$ and this zero is $2^{-1/p}.$ Since%

\begin{align*}
g\left(  0\right)   &  =\allowbreak2,\\
g\left(  1\right)  \allowbreak &  =\allowbreak2,\\
g\left(  2^{-1/p}\right)   &  =2^{\frac{4}{p}},
\end{align*}
and since $2^{\frac{4}{p}}\geq2$ whenever $2<p\leq4$, we finally conclude
that
\[
\max_{a\in\left[  0,1\right]  }\left\{  \left(  2\left\vert \frac{2a^{p}%
-1}{a^{2}+\left(  1-a^{p}\right)  ^{2/p}}\right\vert ^{2}+\left(  2a\left(
1-a^{p}\right)  ^{\frac{1}{p}}\frac{a^{p-2}+\left(  1-a^{p}\right)
^{\frac{p-2}{p}}}{a^{2}+\left(  1-a^{p}\right)  ^{2/p}}\right)  ^{2}\right)
^{\frac{1}{2}}\right\}  =2^{\frac{2}{p}}\text{.}%
\]

If $q>2$, we have $\ell_{2}\subset\ell_{q}$ and
\[
\left\Vert \cdot\right\Vert _{q}\leq\left\Vert \cdot\right\Vert _{2},
\]
and thus%
\begin{align*}
&  C_{2,2,p,q}\\
= &  \max_{a\in\left[  0,1\right]  }\left\{  \left(  2\left\vert \frac
{2a^{p}-1}{a^{2}+\left(  1-a^{p}\right)  ^{2/p}}\right\vert ^{q}+\left(
2a\left(  1-a^{p}\right)  ^{\frac{1}{p}}\frac{a^{p-2}+\left(  1-a^{p}\right)
^{\frac{p-2}{p}}}{a^{2}+\left(  1-a^{p}\right)  ^{2/p}}\right)  ^{q}\right)
^{\frac{1}{q}}\right\}  \\
&  \leq\max_{a\in\left[  0,1\right]  }\left\{  \left(  2\left\vert
\frac{2a^{p}-1}{a^{2}+\left(  1-a^{p}\right)  ^{2/p}}\right\vert ^{2}+\left(
2a\left(  1-a^{p}\right)  ^{\frac{1}{p}}\frac{a^{p-2}+\left(  1-a^{p}\right)
^{\frac{p-2}{p}}}{a^{2}+\left(  1-a^{p}\right)  ^{2/p}}\right)  ^{2}\right)
^{\frac{1}{2}}\right\}  \\
&  =2^{\frac{2}{p}}.
\end{align*}
Now we just need to recall that%
\begin{align*}
&  \left(  2\left\vert \frac{2\left(  2^{-1/p}\right)  ^{p}-1}{a^{2}+\left(
1-\left(  2^{-1/p}\right)  ^{p}\right)  ^{2/p}}\right\vert ^{q}+\left(
2\left(  2^{-1/p}\right)  \left(  1-\left(  2^{-1/p}\right)  ^{p}\right)
^{\frac{1}{p}}\frac{\left(  2^{-1/p}\right)  ^{p-2}+\left(  1-\left(
2^{-1/p}\right)  ^{p}\right)  ^{\frac{p-2}{p}}}{\left(  2^{-1/p}\right)
^{2}+\left(  1-\left(  2^{-1/p}\right)  ^{p}\right)  ^{2/p}}\right)
^{q}\right)  ^{\frac{1}{q}}\\
&  =2^{\frac{2}{p}}%
\end{align*}
to complete the proof.\bigskip

\section{Final comments}

From the previous section we can conclude that if $2<p\leq4$ and $1\leq q<2$,
then $C_{2,2,p,q}>2^{\frac{2}{p}}$. The Theorem \ref{800} seems somewhat
surprising since Araujo \textit{et al}. \cite{aaww} mention in their Remark
4.4 that the \textquotedblleft function cannot be optimized explicitly in
general\textquotedblright. The possibility or not of finding a similar closed
formula for the case $p>4$ is still open and seems to be an interesting
problem. We finish the paper by remarking that our main theorem recovers the
computer assisted numerical table presented in \textit{(\cite{aaww})}, and
also corrects the rounding errors since our estimates are exact. As a very
particular case we conclude that the optimal Hardy--Littlewood constant for
$2$-homogeneous polynomials in $\ell_{2}(\mathbb{R}^{2})$ is $\sqrt{2}$ (this
result was obtained, numerically, as a lower bound in \cite{waaa} and proved
to be sharp in the aforementioned paper of Araujo \textit{et al}.
\textit{\cite{aaww}).}


\begin{thebibliography}{99}                                                                                               %


\bibitem {abps}N. Albuquerque, F. Bayart, D. Pellegrino and J. B.
Seoane-Sep\'{u}lveda, Sharp generalizations of the multilinear
Bohnenblust--Hille inequality, J. Funct. Anal. \textbf{266} (2014), 3726--3740.

\bibitem {aaww}G. Ara\'{u}jo, P. Jim\'{e}nez, G. Mu\~{n}oz-Fern\'{a}ndez, J.
Seoane-Sep\'{u}lveda, Estimates on the norm of polynomials and applications,
Estimates on the norm of polynomials and applications, arXiv:1507.01431v1, 6
Jul 2015.

\bibitem {ara2}G. Ara\'{u}jo, D. Pellegrino, Lower bounds for the constants of
the Hardy-Littlewood inequalities, Linear Algebra Appl. \textbf{463} (2014), 10--15.

\bibitem {apd}G. Ara\'{u}jo, D. Pellegrino and D. da Silva e Silva, On the
upper bounds for the constants of the Hardy-Littlewood, J. Funct. Anal.
\textbf{267} (2014), no. 6, 1878--1888.

%\bibitem {bayart}F. Bayart. \emph{Hardy spaces of Dirichlet series and their
%composition operators}, Monatsh. Math. (3), 136 (2002), 203-236.


\bibitem {bps}F. Bayart, D. Pellegrino and J. B. Seoane-Sep\'{u}lveda, The
Bohr radius of the $n$-dimensional polydisk is equivalent to $\sqrt{(\log
n)/n}$, Adv. Math. \textbf{264} (2014), 726--746.

\bibitem {bh}H. F. Bohnenblust and E. Hille, On the absolute convergence of
Dirichlet series, Ann. of Math. \textbf{32} (1931), 600--622.

\bibitem {mic}G. Botelho, C. Michels, D. Pellegrino, Complex interpolation and
summability properties of multilinear operators. Rev. Mat. Complut. 23 (2010),
no. 1, 139--161.

%\bibitem {diaz}A. Defant, J.C. Diaz, D. Garcia, M. Maestre,
%\emph{Unconditional basis and Gordon-Lewis constants for spaces of
%polynomials}, J. Funct. Anal. 181 (2001), 119--145.

\bibitem {cara}D. Carando, A. Defant and P. Sevilla-Peris, The
Bohnenblust-Hille inequality combined with an inequality of Helson,  to appear
in Proc. Amer. Math. Soc.


\bibitem {waaa}W. Cavalcante, D. N\'{u}\~{n}ez-Alarc\'{o}n, D. Pellegrino, New
lower bounds for the constants in the real polynomial Hardy--Littlewood
inequality, arXiv:1506.00159.



\bibitem {ann}A. Defant, L. Frerick, J. Ortega-Cerda, M. Ouna{\"{\i}}es, K.
Seip, The Bohnenblust-Hille inequality for homogeneous polynomials is
hypercontractive, Ann. of Math. (2), \textbf{174} (2011), 485--497.

\bibitem {def33}A. Defant, P. Sevilla-Peris, A new multilinear insight on
Littlewood's 4/3-inequality, J. Funct. Anal. 256 (2009), no. 5, 1642--1664.

\bibitem {sur}A. Defant, P. Sevilla-Peris, The Bohnenblust-Hille cycle of
ideas from a modern point of view. Funct. Approx. Comment. Math. \textbf{50}
(2014), no. 1, 55--127.

\bibitem {dimant}V. Dimant and P. Sevilla--Peris, Summation of coefficients of
polynomials on $\ell_{p}$ spaces, arXiv:1309.6063v1 [math.FA], to appear in
Publ. Mat.

%\bibitem {dineen}S. Dineen, \emph{Complex analysis on infinite-dimensional
%spaces,} Springer Monographs in Mathematics. Springer-Verlag London, Ltd.,
%London, 1999.


\bibitem {gal}D. Galicer, S. Muro, P. Sevilla-Peris, Asymptotic estimates on
the von Neumann inequality for homogeneous polynomials, arXiv:1504.05547v2, to
appear in Journal fur die reine und angewandte Mathematik (Crelle's Journal)

\bibitem {grecu}B. Grecu, Geometry of 2-homogeneous polynomials on $\ell_{p}$
spaces, $1<p<\infty$. J. Math. Anal. Appl. 273 (2002), no. 2, 262--282.

\bibitem {hardy}G. Hardy and J. E. Littlewood, Bilinear forms bounded in space
$[p,q]$, Quart. J. Math. \textbf{5} (1934), 241--254.

\bibitem {LLL}J.E. Littlewood, On bounded bilinear forms in an infinite number
of variables, Quart. J. (Oxford Ser.) \textbf{1} (1930), 164--174.

\bibitem {monta}A. Montanaro, Some applications of hypercontractive
inequalities in quantum information theory. J. Math. Phys. 53 (2012), no.
\textbf{12}, 122206, 15 pp.

%\bibitem {muno}G. A. Mu\~{n}oz-Fernandez, Y. Sarantopoulos, A. Tonge,
%\emph{Complexifications of real Banach spaces, polynomials and multilinear
%maps}, Studia Math. 134 (1999), 1--33.


\bibitem {pell}D. Pellegrino, The optimal constants of the mixed $\left(
\ell_{1},\ell_{2}\right)  $-Littlewood inequality, to appear in the Journal of
Number Theory, DOI 10.1016/j.jnt.2015.08.007.

\bibitem {popa}D. Popa, G. Sinnamon, Blei's inequality and coordinatewise
multiple summing operators, Publ. Mat. \textbf{57} (2013), no. 2, 455--475.

\bibitem {pra}T. Praciano--Pereira, On bounded multilinear forms on a class of
$\ell_{p}$ spaces, J. Math. Anal. Appl. \textbf{81} (1981), 561--568.

\bibitem {diana}D.M. Serrano-Rodr\'{\i}guez, Improving the closed formula for
subpolynomial constants in the multilinear Bohnenblust--Hille inequalities,
Linear Algebra. Appl. \textbf{438} (2013), 3124--3138.

%\bibitem {weissler}F.B. Weissler, \emph{Logarithmic Sobolev inequalities and¨%hypercontractive estimates on the circle,} J. Funct. Anal. 37 (1980), no. 2, 218--234.

\end{thebibliography}
\end{document}